\documentclass{article}
\newtheorem{theorem}{Theorem}

\newtheorem{definition}[theorem]{Definition}
\newtheorem{example}[theorem]{Example}
\newtheorem{remark}[theorem]{Remark}
\newtheorem{corollary}[theorem]{Corollary}

\def\QED{\quad\blackslug\lower 8.5pt\null}

\newcommand{\crazy}[2]{\displaystyle{\mathop{#1}_{#2}}
\vphantom{\displaystyle{#1}}}

\begin{document}

\setcounter{section}{-1}

\section{Introduction}

Four-dimensional webs $W (3, 2, 2)$ have been considered in many books
and papers (see, for example, the books [AS 92], [G 88]
and the papers [B 35], [C 36], [G 85, 86, 87, 99], [K 81, 83, 84, 96]).
They are of special interest since
\begin{description}
  \item[a)] They are the first webs generalizing the notion of
two-dimensional three-web introduced by Blaschke [Bl 28]
to higher codimension (see [B 35]).

\item[b)] They provide examples illustrating different properties
of webs (see [AS 92],  [B 35], [C 36], [G 88], [G 99]).

  \item[c)] Their torsion tensor has a simple structure: $a^i_{jk}
  = a_{[j} \delta^i_{k]}$, where $a_i$ is a covector (see  [AS 92], [G 88]).

 \item[d)] They are connected with the pseudoconformal structures
  $CO (2, 2)$  of signature $(2, 2)$ (see  [AG 96], [AG 99], [K 81, 83, 96]).
\end{description}

If the covector $a = \{a_1, a_2\}$ of a web  $W (3, 2, 2)$
does not vanish, then it defines a transversal $a$-distribution invariantly
and intrinsically
connected with a web. In general, this $a$-distribution is
not integrable.

In Section 1 we find necessary and sufficient conditions
of its integrability and prove the existence theorem for
 webs $W (3, 2, 2)$ with integrable transversal
$a$-distributions (see Theorems 1 and 3).

In Section 2 we prove that for  a web
$W (3, 2, 2)$ with the  integrable
distribution $\Delta$, its integral surfaces $V^2$  are
geodesicly parallel
in an affine connection of a certain bundle of affine connections (Theorem
4 (i)) and  study three-webs for which the surfaces $V^2$ are
geodesicly parallel with respect to affine connections of this bundle  (Theorem
4 (ii)).

In Section 3, we  find conditions for webs $W (3, 2, 1)$ cut by the foliations
of $W (3, 2, 2)$ on $V^2$ to be hexagonal (Theorem 6) and
prove the existence theorem for such webs  $W (3, 2, 2)$ (Theorem
7). We also  prove the existence theorem for webs $W (3, 2, 2)$
of the subclass  which is
the intersection of subclasses considered in Sections
{\bf 2} and {\bf 3} (Theorem 8),
and establish some properties of webs $W (3, 2, 2)$ implied
by a relationship existing between four-dimensional three-webs
and pseudoconformal structures  $CO (2, 2)$  of signature $(2,
2)$ (Theorem 9).

In addition, in Sections {\bf 2} and {\bf 3} we find an analytic
characterization of three-webs
considered in these sections not only in a specialized frame but
also in the general frame.

Note that
webs $W (3, 2, 2)$ with integrable transversal $a$-distributions
as well as   webs $W (3, 2, 2)$, for
which integral surfaces $V^2$ of $\Delta$  are
geodesicly parallel in an affine connection of a certain bundle of
affine connections, and   webs \linebreak $W (3, 2, 2)$, for which the
three-subwebs  $W (3, 2, 1)$ cut by the foliations
of \linebreak  $W (3, 2, 2)$ on $V^2$ are hexagonal,
are considered in this paper for the first time.

 \section{The transversal distribution of a web $W (3, 2, 2)$}

{\bf 1.} The leaves of the foliation $\lambda_u, \; u =
1, 2, 3$,  of a web $W (3, 2, 2)$ are determined by the equations
$\crazy{\omega}{u}^i = 0, \; i = 1, 2$,   where
\begin{equation}\label{eq:1}
\crazy{\omega}{1}^i + \crazy{\omega}{2}^i + \crazy{\omega}{3}^i = 0
\end{equation}
(see, for example, [G 88], Section {\bf 8.1} or
[AS 92], Section  {\bf 1.3}). The forms $\crazy{\omega}{1}^i$ and
$\crazy{\omega}{2}^i$ are basis forms on a manifold $M^4$ carrying
the web $W (3, 2, 2)$.

The structure equations of such a web can be written in the form
\begin{equation}\label{eq:2}
\renewcommand{\arraystretch}{1.3}
 \left\{
\begin{array}{ll}
   d \crazy{\omega}{1}^i =
\crazy{\omega}{1}^j \wedge \omega_j^i +  a_j \crazy{\omega}{1}^j
\wedge \crazy{\omega}{1}^i, \\  d\crazy{\omega}{2}^i =
\crazy{\omega}{2}^j \wedge \omega_j^i - a_j \crazy{\omega}{2}^j
\wedge \crazy{\omega}{2}^i.
\end{array}
      \right.
\renewcommand{\arraystretch}{1}
 \end{equation}
 The differential prolongations of equations (2) are
(see [G 88], Sections {\bf 8.1} and {\bf 8.4} or [AS 92], Section
{\bf 3.2}):
\begin{equation}\label{eq:3}
  d\omega_j^i - \omega_j^k \wedge
\omega_k^i =   b_{jkl}^i \crazy{\omega}{1}^k \wedge
\crazy{\omega}{2}^l,
 \end{equation}
\begin{equation}\label{eq:4}
da_i - a_j \omega_i^j = p_{ij} \crazy{\omega}{1}^j  +  q_{ij}
\crazy{\omega}{2}^j,
\end{equation}
 where
 \begin{equation}\label{eq:5}
b^i_{[j|l|k]} = \delta^i_{[k} p_{j]l}, \;\; b^i_{[jk]l} =
\delta^i_{[k} q_{j]l}.
\end{equation}
The quantities
\begin{equation}\label{eq:6}
 a_{jk}^i = a_{[j} \delta_{k]}^i
 \end{equation}
 and $b^i_{jkl}$ are the {\em
torsion and curvature tensors} of a three-web $W (3,
2, 2)$. Note that for webs $W (3, 2, 2)$ the torsion tensor
$a^i_{jk}$ always has structure (6), where  $a = \{a_1, a_2\}$ is
its transversal covector. If $a = 0$, then a web $W (3, 2, 2)$ is isoclinicly
geodesic. Such webs were studied in [A 69]. In what follows, {\em we will assume that $a
\neq 0$, i.e., a web $W (3, 2, 2)$ is nonisoclinicly geodesic.}

The covector $a_i$ is defined in a second-order differential
neighborhood of a point $x \in M^4$, and the curvature tensor
$b^i_{jkl}$ as well as the tensors $p_{ij}$ and $q_{ij}$
are  defined in a third-order neighborhood of  $x \in M^4$.
By conditions (5), the tensor $b^i_{jkl}$ can be represented in
the form
$$
b^i_{jkl} = s^i_{jkl} + \frac{2}{3} p_{jk} \delta^i_l
- \frac{1}{3} p_{kl} \delta^i_j
- \frac{1}{3} p_{lj} \delta^i_k
- \frac{1}{3} q_{jk} \delta^i_l
- \frac{1}{3} q_{kl} \delta^i_j
+ \frac{2}{3} q_{lj} \delta^i_k,
$$
where $ s^i_{jkl} = b^i_{(jkl)}$ is the symmetric part
of the tensor $b^i_{jkl}$ (see [AS 92], p. 113). The last formula implies that in
 a third-order neighborhood of  $x \in M^4$, there are 8 independent
 components of the tensors $p_{ij}$ and $q_{ij}$ and also
8 independent  components of the tensor $b^i_{jkl}$.

In this paper we will need the differential prolongations
of equations (3), (4), and (5). They have the form
\begin{equation}\label{eq:7}
  [\nabla b^i_{jkl} + b^i_{jkl} a_m (\crazy{\omega}{1}^m - \crazy{\omega}{2}^m)] \wedge
\crazy{\omega}{1}^k \wedge \crazy{\omega}{2}^l = 0,
 \end{equation}
\begin{equation}\label{eq:8}
  (\nabla p_{jk} + p_{jk} a_l \crazy{\omega}{1}^l) \wedge
\crazy{\omega}{1}^k +
 (\nabla q_{jk} - q_{jk} a_l \crazy{\omega}{2}^l) \wedge
\crazy{\omega}{2}^k + a_m b^m_{jkl}  \crazy{\omega}{1}^k
\wedge  \crazy{\omega}{2}^l = 0,
 \end{equation}
 where
 $$
\renewcommand{\arraystretch}{1.3}
\begin{array}{ll}
\nabla b^i_{jkl} = d b^i_{jkl} - b^i_{mkl}  \omega^m_j
 - b^i_{jml}  \omega^m_k  - b^i_{jkm}  \omega^m_l
 + b^m_{jkl}  \omega^i_m, \\
\nabla p_{jk} = d p_{jk} - p_{mk}  \omega^m_j
 - p_{jm}  \omega^m_k, \\
\nabla q_{jk} = d q_{jk} - q_{mk}  \omega^m_j
 - q_{jm}  \omega^m_k.
\end{array}
\renewcommand{\arraystretch}{1}
 $$

 Equations (7) and (8) prove that the forms
$\nabla b^i_{jkl}, \nabla p_{jk}$, and $\nabla q_{jk}$
are linear combinations of the basis forms
$\crazy{\omega}{1}^k$ and $\crazy{\omega}{2}^k$:
\begin{equation}\label{eq:9}
\renewcommand{\arraystretch}{1.3}
 \left\{
\begin{array}{ll}
  \nabla b^i_{jkl} = \overline{b}^i_{jklm} \crazy{\omega}{1}^m
  +  \widetilde{b}^i_{jklm} \crazy{\omega}{2}^m, \\
  \nabla p_{jk} = \overline{p}_{jkl}  \crazy{\omega}{1}^l
  +  \widetilde{p}_{jkl} \crazy{\omega}{2}^l, \\
 \nabla q_{jk} = \overline{q}_{jkl}  \crazy{\omega}{1}^l
  +  \widetilde{q}_{jkl} \crazy{\omega}{2}^l.
 \end{array}
      \right.
\renewcommand{\arraystretch}{1}
 \end{equation}
Substituting decompositions (9) into equations (7) and (8)
and using the linear independence of the forms
$\crazy{\omega}{\alpha}^i$, we find that the coefficients in (9)
satisfy the following conditions:
\begin{equation}\label{eq:10}
\renewcommand{\arraystretch}{1.3}
 \left\{
\begin{array}{ll}
\overline{b}^i_{j[k|l|m]} + a_{[m} b^i_{|j|k]l} = 0, \\
\widetilde{b}^i_{jk[lm]} - a_{[m} b^i_{|jk|l]} = 0, \\
\overline{p}_{i[lk]} + q_{i[l} a_{k]} = 0, \\
\widetilde{q}_{i[lk]} - q_{i[l} a_{k]} = 0, \\
 a_{m} b^m_{jkl} - \widetilde{p}_{jkl} + \overline{q}_{jlk} = 0.
  \end{array}
      \right.
\renewcommand{\arraystretch}{1}
 \end{equation}
 In addition, upon differentiating conditions (5) and applying
 equations (9), we find other conditions for
 the coefficients in (9):
\begin{equation}\label{eq:11}
\renewcommand{\arraystretch}{1.3}
 \left\{
\begin{array}{ll}
\overline{b}^i_{[j|l|k]m}  = \delta^i_{[k} \overline{p}_{j]lm}, &
\widetilde{b}^i_{[j|l|k]m} = \delta^i_{[k} \widetilde{p}_{j]lm},\\
\overline{b}^i_{[jk]lm}    = \delta^i_{[k}   \overline{q}_{j]lm}, &
\widetilde{b}^i_{[jk]lm}   = \delta^i_{[k} \widetilde{q}_{j]lm}.
  \end{array}
      \right.
\renewcommand{\arraystretch}{1}
 \end{equation}

It follows from conditions (5) and (10)  that
in a fourth-order neighborhood of  $x \in M^4$, there are 6 independent
 components of the tensors $\overline{p}_{ijk}, \widetilde{p}_{ijk},
 \overline{q}_{ijk}, \widetilde{q}_{ijk}$ and also
20 independent  components  of the tensors $\overline{b}^i_{jkl},
 \widetilde{b}^i_{jkl}$.

{\bf 2}. For a web $W (3, 2, 2)$, a transversally geodesic
distribution is defined (cf. [AS 92], Section {\bf 3.1}) by the
equations $$ \xi^2 \crazy{\omega}{1}^1 - \xi^1
\crazy{\omega}{1}^2= 0, \;\; \xi^2 \crazy{\omega}{2}^1 - \xi^1
\crazy{\omega}{2}^2 = 0. $$  If we take $\displaystyle
\frac{\xi^1}{\xi^2} = - \frac{a_2}{a_1}$, we obtain the
invariant transversal distribution $\Delta$ defined by the
equations
\begin{equation}\label{eq:12}
 a_1 \crazy{\omega}{1}^1 + a_2 \crazy{\omega}{1}^2 = 0, \;\; a_1
 \crazy{\omega}{2}^1 + a_2 \crazy{\omega}{2}^2 = 0.
  \end{equation}

This distribution is defined by the  1-forms
\begin{equation}\label{eq:13}
\crazy{\omega}{\alpha} = a_1 \crazy{\omega}{\alpha}^1 + a_2
\crazy{\omega}{\alpha}^2, \;\; \alpha = 1, 2.
\end{equation}
It is  connected with a web invariantly and intrinsically since
it is defined by the torsion tensor of a web.
We will call the distribution $\Delta$ the {\em
transversal $a$-distribution} of a web $W (3, 2, 2)$.
Note that for isoclinicly geodesic webs $W (3, 2, 2)$, for which
$a_1 = a_2 =0$, the  distribution $\Delta$ is not defined.

The following theorem gives the conditions of integrability of the
distribution $\Delta$.

\begin{theorem}
The transversal $a$-distribution $\Delta$ defined by
the equations $(12)$  is integrable if and only if
\begin{equation}\label{eq:14}
\renewcommand{\arraystretch}{1.3}
\left\{
\begin{array}{ll}
a_2^2 p_{11} -  2a_1 a_2 p_{(12)} + a_1^2 p_{22} = 0, \\
a_2^2 q_{11} -  2a_1 a_2 q_{(12)} + a_1^2 q_{22} = 0.
        \end{array}
\right.
\renewcommand{\arraystretch}{1}
\end{equation}
\end{theorem}

{\sf Proof.} A transversal distribution $\Delta$ defined by equations (12)
 is integrable if and only if
\begin{equation}\label{eq:15}
d  \crazy{\omega}{\alpha} \wedge \crazy{\omega}{1} \wedge
\crazy{\omega}{2} = 0, \;\;\;\;\; \alpha = 1, 2.
\end{equation}
By (2) and (13), equations (15) take the forms:
$$
\renewcommand{\arraystretch}{1.3}
\left\{
\begin{array}{ll}
(a_2 \nabla a_1 - a_1  \nabla a_2) \wedge \crazy{\omega}{1}^1
\wedge \crazy{\omega}{1}^2 \wedge (a_1 \crazy{\omega}{2}^1 + a_2
\crazy{\omega}{2}^2) = 0, \\ (a_2 \nabla a_1 - a_1  \nabla a_2)
\wedge \crazy{\omega}{2}^1 \wedge \crazy{\omega}{2}^2 \wedge (a_1
\crazy{\omega}{1}^1 + a_2 \crazy{\omega}{1}^2) = 0,
        \end{array}
\right.
\renewcommand{\arraystretch}{1}
$$
where $\nabla a_i = d a_i - a_j \omega^j_i$.
By (4) and the linear independence of the forms
$\crazy{\omega}{\alpha}^i$,   the last equations
imply conditions (14).
 \rule{3mm}{3mm}

Note that for an arbitrary web $W (3, 2, 2)$, it is always possible
to take a specialized frame in which  there is a relation
between the components $a_1$ and $a_2$  of the covector $a$.
For example, if the transversal distribution $\Delta$ coincides
with the distribution $\crazy{\omega}{\alpha}^1 = 0$ or $\crazy{\omega}{\alpha}^2 = 0$
or $\crazy{\omega}{\alpha}^1 + \crazy{\omega}{\alpha}^2 = 0$
or $\crazy{\omega}{\alpha}^1 - \crazy{\omega}{\alpha}^2 = 0$, then we
have $a_2 = 0$ or $a_1 = 0$ or $a_1 = a_2$ or $a_1 = - a_2$, respectively.
Note that in these cases
the forms $\omega_2^1, \; \omega_1^2, \; \omega_1^1 + \omega_1^2
- \omega_2^1 - \omega_2^2$, and $- \omega_1^1 + \omega_1^2
- \omega_2^1 + \omega_2^2$, respectively, are expressed in terms of
the basis forms $\crazy{\omega}{\alpha}^i$, i.e., in these cases we have
$$
\pi_2^1 = 0, \;\; \pi_1^2 = 0, \;\; \pi_1^1 + \pi_1^2 - \pi_2^1 - \pi_2^2 = 0,
\;\; - \pi_1^1 + \pi_1^2 - \pi_2^1 + \pi_2^2 = 0,
$$
respectively, where $\pi_i^j = \omega_i^j\Bigl|_{\crazy{\omega}{\alpha}^i = 0}$.

In proving the existence theorems it is convenient to use
one of these specializations. Let us reformulate Theorem 1 for
the first specialization indicated above.

 \begin{corollary}
 If the frame bundle associated with a three-web
 $W (3, 2, 2)$ is specialized in such a way that
\begin{equation}\label{eq:16}
a_2 = 0,
\end{equation}
then the $a$-distribution $\Delta$ coincides with
the coordinate distribution $\crazy{\omega}{\alpha}^1 = 0$ and
the condition $ \pi_2^1 = 0$ holds. In such a frame bundle
the $a$-distribution is integrable
if and only if
\begin{equation}\label{eq:17}
p_{22} = 0, \;\; q_{22} = 0.
\end{equation}

\end{corollary}
{\sf Proof}. This follows from equations (14) and (16).
\rule{3mm}{3mm}

Each of the relations (14) and (17) gives two
conditions which Pfaffian derivatives $p_{ij}$ and $q_{ij}$ of the covector $a$ must
satisfy in order for the transversal distribution $\Delta$
of a web $W (3, 2, 2)$ to be integrable.

\vspace*{3mm}

{\bf 3.} We will now prove an  existence theorem
for webs with  integrable transversal $a$-distributions $\Delta$.

\begin{theorem} The webs with  integrable  transversal $a$-distributions $\Delta$
exist, and a solution of a system of differential equations
defining such webs depends on five arbitrary functions of three
variables.
\end{theorem}

{\sf Proof.} Suppose that specialization (16) has been made.
Since our web $W (3, 2, 2)$ has the integrable $a$-distribution
$\Delta$ defined by the equations $\crazy{\omega}{\alpha}^1 = 0$,
we have conditions  (16) and (17), and
equations (4) take the form
\begin{equation}\label{eq:18}
\renewcommand{\arraystretch}{1.3}
\left\{
      \begin{array}{ll}
da_1 - a_1 \omega_1^1 = p_{1j} \crazy{\omega}{1}^j + q_{1j} \crazy{\omega}{2}^j,
 \\    - a_1 \omega_2^1 = p_{21} \crazy{\omega}{1}^1 + q_{21}
     \crazy{\omega}{2}^1.
     \end{array}
     \right.
\renewcommand{\arraystretch}{1}
\end{equation}

The exterior cubic and quadratic  equations (7) and (8) become
\begin{equation}\label{eq:19}
  [\nabla b^i_{jkl} + b^i_{jkl} a_1 (\crazy{\omega}{1}^1 - \crazy{\omega}{2}^1] \wedge
\crazy{\omega}{1}^k \wedge \crazy{\omega}{2}^l = 0,
 \end{equation}
\begin{equation}\label{eq:20}
\renewcommand{\arraystretch}{1.3}
\left\{
      \begin{array}{ll}
  (\nabla p_{1k} + p_{1k} a_1 \crazy{\omega}{1}^1) \wedge
\crazy{\omega}{1}^k + (\nabla q_{1k} - q_{1k} a_1 \crazy{\omega}{2}^1) \wedge
\crazy{\omega}{2}^k + a_1 b^1_{1kl}  \crazy{\omega}{1}^k
\wedge  \crazy{\omega}{2}^l = 0,\\
   \nabla p_{21} \wedge \crazy{\omega}{1}^1
   + \nabla q_{21} \wedge
\crazy{\omega}{2}^1 + a_1 b^1_{2kl}  \crazy{\omega}{1}^k
\wedge  \crazy{\omega}{2}^l = 0.
       \end{array}
      \right.\renewcommand{\arraystretch}{1.3}
 \end{equation}

First note that the last equation of (20) implies that
\begin{equation}\label{eq:21}
 b^1_{222} = 0.
  \end{equation}

By (17), (21), and (5), the number of unknown 1-forms (6 forms
$\nabla p_{1i}, \nabla q_{1i}, \linebreak  \nabla p_{21},
\nabla q_{21}$ and   7  forms
among the forms $\nabla b^i_{jkl}$)
is 13,  $q = 13$ (see [BCGGG 91]).

Since we have 2 exterior quadratic equations and 4 exterior
cubic equations (see (19) and (20)), the Cartan's characters are: $s_1 = 2,
s_2 = 6$, and $s_3 = 13 - 8 = 5$. As a result, we have
$Q = s_1 + 2 s_2 + 3 s_3 = 29$.

By (10), 13 Pfaffian derivatives of the functions $p_{1i}, q_{1i},
p_{21}$ and $q_{21}$ are independent: 3 functions
$\overline{p}_{1jk}$, 4 functions $\widetilde{p}_{1jk}$,
 3 functions $\widetilde{q}_{1jk}$, and 3 functions $\overline{p}_{211}, \widetilde{p}_{211},
 \widetilde{q}_{211}$.
In addition, by (5), (10), and (11),  there are 16 independen\mbox{t} functions
 $\overline{b}^i_{jklm}$ and $\widetilde{b}^i_{jklm}$: 10 functions among
 $\overline{b}^2_{jklm}, \widetilde{b}^2_{jklm}$ and 6 function\mbox{s}
 $\overline{b}^1_{1111}, \overline{b}^1_{1112}, \overline{b}^1_{1122},
  \widetilde{b}^1_{1111}, \widetilde{b}^1_{1112},
  \widetilde{b}^1_{1122}$. This implies that
  the general third-order integral element depends on  $N = 13 + 16 = 29$
  parameters.

  Thus, we have $Q = N$. As a result, the system defining three-webs
$W (3, 2, 2)$ is in involution, and its solution depends on
five arbitrary functions of three variables  (see [BCGGG 91]). \rule{3mm}{3mm}

 \section{Geodesicity of integral surfaces}

{\bf 1.} Suppose that the specialization of frames indicated in
Section {\bf 1} has been made, i.e., we have condition (16):
$a_2 = 0$. Then by (12), the distribution $\Delta$ is determined by the
system of equations

\begin{equation}\label{eq:22}
\crazy{\omega}{1}^1 = 0, \;\; \crazy{\omega}{2}^1 = 0.
\end{equation}

In $T_x (M)$, consider a vectorial frame \{$e_i^\alpha$\} that is conjugate to
the coframe $\{\crazy{\omega}{\alpha}^i\}$. Thus for $x \in M$, we
obtain
$$
d x = e_i^\alpha \crazy{\omega}{\alpha}^i.
$$
Then on integral surfaces $V^2$ of the $a$-distribution $\Delta$,
we have
\begin{equation}\label{eq:23}
d x = e_2^1 \crazy{\omega}{1}^2 + e_2^2 \crazy{\omega}{2}^2.
\end{equation}
The 1-forms $\crazy{\omega}{1}^2$ and $\crazy{\omega}{2}^2$ are
basis forms on surfaces $V^2$, and the vectors $e_2^1$ and $e_2^2$
 are tangent to these surfaces.

Consider the affine connections $\Gamma$ defined on the manifold $M^4$ by
1-forms
\begin{equation}\label{eq:24}
  \theta_v^u = \pmatrix{\theta_j^i & 0 \cr
               0 & \theta_j^i \cr}, \;\;\;\; i, j = 1, 2;
               \;\; u, v = 1, 2, 3, 4,
\end{equation}
where
\begin{equation}\label{eq:25}
  \theta_j^i = \omega_j^i + a^i_{jk} (p \crazy{\omega}{1}^k
  + q\crazy{\omega}{2}^k)
\end{equation}
(see [AS 92], p. 35). For the three-web $W (3, 2, 2)$ in question,
by (6) and (16), formulas (25) take the form:
\begin{equation}\label{eq:26}
\renewcommand{\arraystretch}{1.3}
\left\{
\begin{array}{ll}
  \theta_1^1 = \omega_1^1, &  \theta_1^2 = \omega_1^2 +
\frac{1}{2} a_1 (p \crazy{\omega}{1}^2 + q\crazy{\omega}{2}^2), \\
 \theta_2^1 = \omega_2^1, &  \theta_2^2 = \omega_2^2 +
\frac{1}{2} a_1 (p \crazy{\omega}{1}^1 + q\crazy{\omega}{2}^1).
 \end{array}
\right.
\renewcommand{\arraystretch}{1}
\end{equation}

When the vectorial frame \{$e_i^\alpha$\} moves along the manifold $M^4$ endowed
with a connection $\Gamma$, we obtain
\begin{equation}\label{eq:27}
\renewcommand{\arraystretch}{1.3}
\left\{
\begin{array}{ll}
d e_2^1 = \theta_2^1 e_1^1 + \theta_2^2 e_2^1, \\
d e_2^2 = \theta_2^1 e_1^2 + \theta_2^2 e_2^2.
 \end{array}
\right.
\renewcommand{\arraystretch}{1}
\end{equation}

We will now prove the following result.

\begin{theorem}
\begin{description}
\item[(i)] If the $a$-distribution $\Delta$ is integrable on the
web $W (3, 2, 2)$, then its integral surfaces $V^2$ are totally geodesic
on $M^4$ in any affine connection of the
bundle $(24)$--$(25)$.

\item[(ii)]
If on a web $W (3, 2, 2)$ the  condition
\begin{equation}\label{eq:28}
\omega_2^1 = 0
\end{equation}
holds, then the integral surfaces $V^2$ of the $a$-distribution
$\Delta$ are geodesicly parallel in any affine connection of the
bundle $(24)$--$(25)$.
\end{description}
\end{theorem}

{\sf Proof}.
\begin{description}
\item[(i)] By the second of relations (18) and (22), on surfaces $V^2$ we have
\begin{equation}\label{eq:29}
\omega_2^1|_{V^2} = 0.
\end{equation}
This and equations (22) and (26) imply that on $V^2$ equations
(27) take the form
\begin{equation}\label{eq:30}
d e_2^1 =  \theta_2^2 e_2^1, \;\; d e_2^2 = \theta_2^2 e_2^2.
\end{equation}
It follows that the bivectors $\Delta = e_2^1 \wedge e_2^2$
are geodesicly parallel on $V^2$
in any affine connection of the
bundle (24)--(25). As a result, the integral surfaces $V^2$ are totally geodesic
on $M^4$ in any of these connections.

\item[(ii)] If equations (28) hold on the entire manifold $M^4$,
then equations (30) are identically satisfied on $M^4$. Therefore,
the bivectors $\Delta = e_2^1 \wedge e_2^2$ are geodesicly parallel on the
entire manifold $M^4$ in any affine connection of the
bundle $(24)$-$(25)$. As a result, the integral surfaces $V^2$ of the
$a$-distribution $\Delta$ are not only totally geodesic but also
 geodesicly parallel on the
entire manifold $M^4$ in all these connections.
\end{description}

{\bf 2.} Preliminary considerations show that
 three-webs $W (3, 2, 2)$, for which integral surfaces $V^2$
of the transversal $a$-distribution $\Delta$ are  geodesicly
parallel, exist, and a solution of a system defining such webs
depends on four arbitrary functions of three variables. However,
we were not able to check the Cartan test in detail.

{\bf 3.}  Conditions (28) for
 integral surfaces $V^2$  of the transversal $a$-distribution
 $\Delta$ to be geodesicly parallel were obtained in a specialized
 frame, i.e., for $a_2 = 0$. To find these conditions in the
 general frame, we first note that by (18), equations (28) are
 equivalent to equations
\begin{equation}\label{eq:31}
p_{21} = 0, \;\; q_{21} = 0
\end{equation}
Of course, conditions (17) of integrability of the $a$-distribution
$\Delta$ in a specialized frame must be added to conditions (31).

In order to write equations (17) and (31) in the
 general frame, we  write equations (13) in the form
\begin{equation}\label{eq:32}
 \crazy{\omega}{\alpha}^{1'} = a_1  \crazy{\omega}{\alpha}^{1} +
 a_2  \crazy{\omega}{\alpha}^{2},
\end{equation}
consider a relation
\begin{equation}\label{eq:33}
 \crazy{\omega}{\alpha}^{2'} = c_1  \crazy{\omega}{\alpha}^{1} +
 c_2  \crazy{\omega}{\alpha}^{2},
\end{equation}
along with equation (32), and  assume that
\begin{equation}\label{eq:34}
D = \det \pmatrix{a_1 & a_2 \cr
                  c_1  & c_2 \cr} \neq 0.
\end{equation}
The 1-forms $ \crazy{\omega}{\alpha}^{1'}$ and
$\crazy{\omega}{\alpha}^{2'}$ form a basis on a manifold $M^4$
carrying a three-web $W (3, 2, 2)$ whose coordinate bivectors
determined by the equations $ \crazy{\omega}{\alpha}^{1'} = 0$ and
$\crazy{\omega}{\alpha}^{2'} = 0$ are transversal. The first
of these bivectors is defined by the torsion tensor of
the three-web $W (3, 2, 2)$, and the second one is chosen
arbitrarily.

Let us write equations (32) and (33) in the form
\begin{equation}\label{eq:35}
 \crazy{\omega}{\alpha}^{i'} = a_j^{i'}  \crazy{\omega}{\alpha}^{j},
 \end{equation}
where the matrix
\begin{equation}\label{eq:36}
A = (a_j^{i'}) =  \pmatrix{a_1 & a_2 \cr
                  c_1  & c_2\cr}
\end{equation}
 is nondegenerate. Its inverse matrix can be written in the form
\begin{equation}\label{eq:37}
A^{-1} = (a_{i'}^j) =  \frac{1}{D}  \pmatrix {c_2 & -a_2 \cr
                  - c_1  & a_1\cr}.
 \end{equation}

 Under the coframe transformation (35),
the tensors $p_{ij}$ and $q_{ij}$ of the web $W (3, 2, 2)$
undergo the regular  tensor  transformation:
\begin{equation}\label{eq:38}
p_{i'j'} =  a_{i'}^i a_{j'}^j p_{ij}, \;\;
q_{i'j'} =  a_{i'}^i a_{j'}^j q_{ij}.
 \end{equation}
Taking into account (37), we write formulas (38)
for the components $p_{21}, q_{21},  p_{22}$ and $q_{22}$
of these tensors:
\begin{equation}\label{eq:39}
\renewcommand{\arraystretch}{1.3}
\left\{
\begin{array}{ll}
p_{2' 1'} = c_1 (a_{2} p_{12} - a_1 p_{22}) + c_2 (a_1 p_{21} - a_{2} p_{11}),
\\
q_{2' 1'} = c_1 (a_{2} q_{12} - a_1 q_{22}) + c_2 (a_1 q_{21} - a_{2} q_{11}),\\
p_{2' 2'} = - a_1 (a_{2} p_{12} - a_1 p_{22}) - a_2 (a_1 p_{21} - a_{2} p_{11}),
\\
q_{2' 2'} = - a_1 (a_{2} q_{12} - a_1 q_{22}) - a_2 (a_1 q_{21} - a_{2} q_{11}).
\end{array}
\right.
\renewcommand{\arraystretch}{1}
\end{equation}
Note that the right-hand sides of
the last two expressions differ from the left-hand sides
of equations (14) only by sign.

Conditions (39) imply the following result.

\begin{theorem} The integral surfaces
$V^2$ of the $a$-distribution $\Delta$ are geodesicly parallel
with respect to any affine connection of the bundle $(24)$--$(25)$ if and only
if the components of the covector $a = \{a_i\}$ and of the tensors
$p_{ij}$ and $q_{ij}$ satisfy the following conditions:
\begin{equation}\label{eq:40}
\renewcommand{\arraystretch}{1.3}
\left\{
\begin{array}{ll}
a_{2} p_{12} - a_1 p_{22} = 0, & a_1 p_{21} -  a_{2} p_{11} = 0, \\
a_{2} q_{12} - a_1 q_{22} = 0, &  a_1 q_{21} - a_{2} q_{11} = 0.
\end{array}
\right.
\renewcommand{\arraystretch}{1}
\end{equation}
\end{theorem}

{\sf Proof.} In fact, by (17) and (31),  necessary and sufficient conditions for
 integral surfaces $V^2$  of the transversal $a$-distribution
 $\Delta$ to be geodesicly parallel in the general frame with
 respect to any affine connection of the bundle (24)--(25) have the form
\begin{equation}\label{eq:41}
p_{2' 1'} = 0, \;\; q_{2' 1'} = 0, \;\;
p_{2' 2'} = 0, \;\; q_{2' 2'} = 0.
\end{equation}
But by conditions (34) and (39), equations (41) are equivalent to conditions
(40). \rule{3mm}{3mm}

\section{Hexagonality of two-dimensional three-subwebs}

{\bf 1}. On integral surfaces of the $a$-distribution
$\Delta$ defined on $M^4$ by the torsion tensor
of a web $W (3, 2, 2)$, the leaves of this web cut
two-dimensional three-subwebs $W (3, 2, 1)$.
Let us find the structure equations of these subwebs.

In a specialized frame in which condition (16) holds, the integral
surfaces $V^2$ are defined by the system of equations (22). In
addition, the complete integrability of this system on a surface
$V^2$ and equations (31) imply that equation (28) holds. Thus,
on a surface $V^2$, we have
\begin{equation}\label{eq:42}
\crazy{\omega}{1}^1 = 0, \;\; \crazy{\omega}{2}^1 = 0, \;\;
\omega_2^1= 0,
 \end{equation}
 and the forms $\crazy{\omega}{1}^2$ and $\crazy{\omega}{2}^2$ are
 basis forms on $V^2$.
One-dimensional foliations of a web $W (3, 2, 1)$ are defined on
$V^2$ by the equations
\begin{equation}\label{eq:43}
\crazy{\omega}{1}^2= 0, \;\; \crazy{\omega}{2}^2 = 0, \;\;
\crazy{\omega}{1}^2 + \crazy{\omega}{2}^2 = 0.
 \end{equation}

 To find the structure equations of  webs $W (3, 2, 1)$ on
 surfaces $V^2$, we substitute the values (31) of the forms
$\crazy{\omega}{1}^1, \,\crazy{\omega}{2}^1 $, and $\omega_2^1$
into equations (2) and (3). As a result, we obtain the following
structure equations:
\begin{equation}\label{eq:44}
\renewcommand{\arraystretch}{1.3}
\left\{
\begin{array}{ll}
d \crazy{\omega}{1}^2 = \crazy{\omega}{1}^2 \wedge \omega_2^2, \\
d \crazy{\omega}{2}^2 = \crazy{\omega}{2}^2 \wedge \omega_2^2, \\
d \omega_2^2 = b^2_{222} \crazy{\omega}{1}^2  \wedge
\crazy{\omega}{2}^2.
\end{array}
\right.
\renewcommand{\arraystretch}{1}
\end{equation}
Comparing these equations with the structure equations of
a two-dimensional three-web (see [AS 92], p. 18), we see that
the form $\omega_2^2$ is the connection form of the web
$W (3, 2, 1)$, and the component $b^2_{222}$ of the curvature tensor
of the web $W (3, 2, 2)$ is the curvature of the web $W (3, 2, 1)$:
$$
K = b^2_{222}.
$$

Since the vanishing of the curvature of the web  $W (3, 2, 1)$ is
a necessary and sufficient condition for its hexagonality, we
arrive at the following result.

\begin{theorem}
Two-dimensional three-webs $W (3, 2, 1)$ cut on integral surfaces
of the $a$-distribution $\Delta$ by the foliations of $W (3, 2, 2)$
are hexagonal if and only if in the specialized frame bundle defined by
condition $(16)$,  the component $b^2_{222}$ of the curvature tensor
of the web $W (3, 2, 2)$ vanishes.
\end{theorem}

{\bf 2.} We will now prove an  existence theorem
for webs, for which   two-dimensional  three-webs $W (3, 2, 1)$ cut
on integral surfaces $V^2$  by the foliations of the web $W (3, 2, 2)$ are
hexagonal.

\begin{theorem} The webs  $W (3, 2, 2)$, for which
  two-dimensional  three-webs $W (3, 2, 1)$ cut
on  integral surfaces $V^2$  by the foliations of the web $W (3, 2, 2)$ are
hexagonal, exist, and a solution of a system of differential equations
defining such webs depends on four arbitrary functions of three
variables.
\end{theorem}

{\sf Proof.} Suppose that  specialization (16) has been made.
Then $a_2 = 0$.  Since the $a$-distribution $\Delta$ is integrable,
we have conditions (17). As a result, equations (4) take the form
(18). As we showed in the proof of Theorem 3,
the second of equations (18) implies (21).

Finally, since two-dimensional  three-webs $W (3, 2, 1)$ cut
on integral surfaces $V^2$  by the foliations of the web $W (3, 2, 2)$, are
hexagonal, we have
\begin{equation}\label{eq:45}
  b^2_{222} = 0.
\end{equation}

By (17), (21), (45), and (5),
there are 4 exterior cubic   equations (7) and
two exterior quadratic equation (8).

By  (17), (21), and  (5), the number of unknown 1-forms (6 forms
$\nabla p_{1i}, \nabla q_{1i},  \linebreak \nabla p_{21}, \nabla q_{21}$
 and 6  forms $\nabla b^i_{jkl}$, namely, the forms
$\nabla b^1_{111}, \nabla b^1_{112}, \nabla b^1_{122}, \nabla
b^2_{111},  \linebreak \nabla b^2_{112}, \nabla b^2_{122}$)
is 12,  $q = 18$ (see [BCGGG 91]).

Thus, the Cartan's characters are: $s_1 = 2,
s_2 = 6$, and $s_3 = 12 - 8 = 4$. As a result, we have
$Q = s_1 + 2 s_2 + 3 s_3 = 26$.

By (10) and (45), 14 Pfaffian derivatives of the
functions $p_{1i}$ and $q_{1i}$
are independent: 4 functions
$\overline{p}_{111}, \overline{p}_{112},\overline{p}_{122},\overline{p}_{211}$,
4 functions $\widetilde{p}_{1jk}$,
and  4 function\mbox{s} $\widetilde{q}_{111}, \widetilde{q}_{112},
 \widetilde{q}_{122},
\widetilde{q}_{211}$.
In addition, by (5), (10), (11),  and (45),  there are 12 independent
functions among $\overline{b}^i_{jklm}$ and $\widetilde{b}^i_{jklm}$:
$\overline{b}^i_{1111}, \overline{b}^i_{1112}, \overline{b}^i_{1122}$,
$\widetilde{b}^i_{1111}, \widetilde{b}^i_{1112}, \widetilde{b}^i_{1122}$.
This implies that
the general third-order integral element depends on $N = 14 + 12 = 26$
parameters.

  Thus, we have $Q = N$. As a result, the system defining
  three-webs, for which
and  two-dimensional  three-webs $W (3, 2, 1)$ cut
on integral surfaces $V^2$  by the foliations of the web $W (3, 2, 2)$ are
hexagonal,  is  in involution, and
its solution depends on
four arbitrary functions of three variables
(see [BCGGG 91]). \rule{3mm}{3mm}

{\bf 3.} We will now prove an  existence theorem
for webs, for which integral surfaces $V^2$
of the transversal distribution $\Delta$ are  geodesicly parallel
and  two-dimensional  three-webs $W (3, 2, 1)$ cut
on $V^2$  by the foliations of the web $W (3, 2, 2)$, are
hexagonal.

\begin{theorem} The webs  $W (3, 2, 2)$, for which integral surfaces $V^2$
of the transversal distribution $\Delta$ are  geodesicly parallel,
and  two-dimensional  three-webs  \linebreak $W (3, 2, 1)$ cut
on $V^2$  by the foliations of the web $W (3, 2, 2)$, are
hexagonal, exist, and a solution of a system of differential equations
defining such webs depends on three arbitrary functions of three
variables.
\end{theorem}

{\sf Proof.} Suppose that  specialization (16) has been made.
Then $a_2 = 0$.  Since the surfaces $V^2$ are  geodesicly parallel, we have
conditions (28), i.e., we have $\omega_2^1 = 0$. As a result, equations (4) take the form
\begin{equation}\label{eq:46}
\renewcommand{\arraystretch}{1.3}
\left\{
      \begin{array}{ll}
da_1 - a_1 \omega_1^1 = p_{1j} \crazy{\omega}{1}^j + q_{1j} \crazy{\omega}{2}^j,
  \\   \omega_2^1 = 0.
     \end{array}
     \right.
\renewcommand{\arraystretch}{1}
\end{equation}
 By  (4), the second of equations (46) implies  that
\begin{equation}\label{eq:47}
  p_{2i} = 0, \;\; q_{2i} = 0,
\end{equation}
and by (3), the same equation implies that
\begin{equation}\label{eq:48}
  b^1_{2kl} = 0.
\end{equation}

Since two-dimensional  three-webs $W (3, 2, 1)$ cut
on $V^2$  by the foliations of the web $W (3, 2, 2)$, are
hexagonal, we have condition (45):
$$
  b^2_{222} = 0.
$$
Note that conditions  (47) imply conditions (17) of integrability
of the distribution $\Delta$ defined by the equations
$\crazy{\omega}{\alpha}^1 = 0$.

By (47), (48), (45),  and (5),
there are 3 exterior cubic   equations (7) and
only one exterior quadratic equation (8).

By  (47), (48), and  (5), the number of unknown 1-forms (4 forms
$\nabla p_{1i}, \nabla q_{1i}$
 and 4  forms $\nabla b^i_{jkl}$, namely, the forms
$\nabla b^1_{111}, \nabla b^2_{111}, \nabla b^2_{112}, \nabla b^2_{122}$)
is 8,  $q = 8$ (see [BCGGG 91]).

Thus, the Cartan's characters are: $s_1 = 1,
s_2 = 4$, and $s_3 = 8 - 5 = 3$. As a result, we have
$Q = s_1 + 2 s_2 + 3 s_3 = 18$.

By (10) and (45), 10 Pfaffian derivatives of the functions $p_{1i}$ and $q_{1i}$
are independent: 3 functions
$\overline{p}_{111}, \overline{p}_{112}, \overline{p}_{122}$,
4 functions $\widetilde{p}_{1jk}$,
and  3 functions $\widetilde{q}_{111}, \widetilde{q}_{112},  \linebreak \widetilde{q}_{122}$,
In addition, by (5), (10), (11), (45), and (48),  there are 8 independent
functions among $\overline{b}^i_{jklm}$ and $\widetilde{b}^i_{jklm}$:
$\overline{b}^1_{1111}, \overline{b}^2_{1111},
\overline{b}^2_{1112}, \overline{b}^2_{1122},
\overline{b}^2_{1222}, \overline{b}^2_{2222},
\widetilde{b}^1_{1111}, \widetilde{b}^2_{1111}$.  \linebreak This implies that
the general third-order integral element depends on $N = 10 + 8 = 18$
parameters.

  Thus, we have $Q = N$. As a result, the system defining
  three-webs, for which integral surfaces $V^2$
of the transversal $a$-distribution $\Delta$ are  geodesicly
parallel,
and  two-dimensional  three-webs $W (3, 2, 1)$ cut
on $V^2$  by the foliations of the web $W (3, 2, 2)$ are
hexagonal,  is  in involution, and
its solution depends on
three arbitrary functions of three variables
(see [BCGGG 91]). \rule{3mm}{3mm}

{\bf 4.} Theorem 6 does not give a condition for
two-dimensional  three-webs $W (3, 2, 1)$ cut
on integral surfaces $V^2$ of the transversal $a$-distribution $\Delta$
 by the foliations of the web $W (3, 2, 2)$ to be
hexagonal in the general frame. To find such a condition in
the general frame,  we note that under the coframe transformation (34),
the curvature tensor of the web $W (3, 2, 2)$ undergoes the regular
 tensor transformation:
\begin{equation}\label{eq:49}
b^{i'}_{j' k' l'} =  a_{i}^{i'} a_{j'}^j a_{k'}^k a_{l'}^l
b^i_{jkl}.
 \end{equation}
We write formulas (49) for the components $b^{1'}_{2' 2' 2'} $ and
$b^{2'}_{2' 2' 2'}$ of the curvature tensor:
\begin{equation}\label{eq:50}
b^{1'}_{2' 2' 2'} = a_{i}^{1'} b^i, \;\;
b^{2'}_{2' 2' 2'} = a_{i}^{2'} b^i,
\end{equation}
where we denote by $b^i$ the following contraction:
$$
b^i = b^i_{jkl} a_{2'}^{j} a^{k}_{2'}  a^{l}_{2'}.
$$
By (37), this contraction can be written as
\begin{equation}\label{eq:51}
b^i = \frac{1}{D^3} (- b^i_{111} a_2^3 + 3  b^i_{(112)} a_2^2 a_1
-  3  b^i_{(122)} a_2 a_1^2 +  b^i_{222} a_1^3).
\end{equation}
Equation (51) shows that the contraction $b^i$ is expressed only
in terms of components of the torsion and curvature tensors of the
web $W (3, 2, 2)$; that is, $b^i$ is completely determined by this
web.

We will now prove the following result.

\begin{theorem} Let $W (3, 2, 2)$ be a four-dimensional three-web
with a nonvanishing covector $a$ and with the integrable
transversal $a$-distributions $\Delta$. Two-dimensional
three-webs $W (3, 2, 1)$ cut
on integral surfaces $V^2$  of the transversal $a$-distribution $\Delta$
 by the foliations of the web $W (3, 2, 2)$ are hexagonal if and only if
the torsion and curvature tensors of this web are connected by the
relations
\begin{equation}\label{eq:52}
b^1 = 0, \;\; b^2 = 0.
\end{equation}
\end{theorem}

{\sf Proof.} Since the $a$-distribution $\Delta$ is integrable, then
in the specialized frame condition  (21) holds.
The hexagonality of the webs $W (3, 2, 1)$
implies that in the specialized frame $b^2_{222} = 0$.
In the general frame these two conditions have the form
\begin{equation}\label{eq:53}
b^{1'}_{2' 2' 2'} = 0, \;\; b^{2'}_{2' 2' 2'} = 0.
\end{equation}
By (36), (37), and (50), equations (53) can be written as follows:
\begin{equation}\label{eq:54}
a_1 b^1 + a_2 b^2 = 0, \;\; c_1 b^1 + c_2 b^2 = 0 .
\end{equation}
Since by (34) $D\neq 0$, equations (54) imply conditions
(53). \rule{3mm}{3mm}

{\bf 5.}
A three-web $W (3, 2, 2)$ defines on a manifold $M^4$ a conformal
structure $CO (2, 2)$ whose isotropic cones $C_x$ in the tangent
space $T_x (M^4)$ are determined by the equation
$$
 \crazy{\omega}{1}^{1} \crazy{\omega}{2}^{2} -
 \crazy{\omega}{1}^{2} \crazy{\omega}{2}^{1} = 0
 $$
 (see [AG 96], p. 196). Transversal bivectors of the three-web
 $W (3, 2, 2)$ form one of two families of planar generators of
 the cones $C_x$. These bivectors are defined by  equations
(12). They can be written in the form
$$
 \crazy{\omega}{1}^{1} + t \crazy{\omega}{1}^{2} = 0, \;\;
 \crazy{\omega}{2}^{1} + t \crazy{\omega}{2}^{2} = 0,
$$
where $t = \displaystyle \frac{a_2}{a_1}$. On the manifold $M^4$, these
bivectors form a fiber bundle $E_\alpha$ whose base is $M^4$
  and whose one-dimensional fibers are defined by the fiber
  parameter $t$.

 The relative conformal curvature of
these bivectors is defined by the formula
\begin{equation}\label{eq:55}
C (t) =  s^2_{111} t^4 -  (3  s^2_{112} -  s^1_{111}) t^3
+ 3 (s^2_{122} -  3 s^1_{112}) t^2  - (3  s^2_{222} - 3 s^1_{122}) t
- s^1_{222},
\end{equation}
where $s^i_{jkl} = b^i_{(jkl)}$ is the symmetrized  curvature
tensor of the web in question ([AG 96], Ch. 5; see also [K 83, 84, 96]).
The vanishing of the
quantity $C (t)$ singles out four transversal bivectors on the cone
$C_x$. These bivectors are called {\em principal}.

Next, on a web $W (3, 2, 2)$ we consider the following
contraction:
\begin{equation}\label{eq:56}
b = b^i a_i.
\end{equation}
This quantity is an absolute invariant of a  web $W (3, 2, 2)$.
Substituting the values (51) of the quantities $b^i$ into
equations (56), we find that
\begin{equation}\label{eq:57}
\renewcommand{\arraystretch}{1.3}
\begin{array}{ll}
b = &\!\!\!\! - \frac{1}{D^3} \bigl[b^2_{111} a_2^4 -  (3  b^2_{(112)} -  b^1_{111})
a_2^3a_1  + 3 (b^2_{(122)} -  b^1_{(112)}) a_2^2 a_1^2 \\
& \!\!\!\! \\
& \!\!\!\! - (b^2_{222} - 3 b^1_{122}) a_2 a_1^3 + b^1_{222} a_1^4\bigr].
     \end{array}
\renewcommand{\arraystretch}{1}
\end{equation}
Comparing equations (55) and (57), we easily find that
\begin{equation}\label{eq:58}
b = - \frac{a_1^4}{D^3} C \Bigl(\frac{a_2}{a_1}\Bigr).
\end{equation}
This means that the invariant $b$ of a  web $W (3, 2, 2)$ differs
 from the relative conformal curvature of the
transversal bivector $\Delta$ defined by the torsion
tensor of $W (3, 2, 2)$ only by a factor.

Relations (50) and (56) allow us to prove the following result.

\begin{theorem} Let $W (3, 2, 2)$ be a four-dimensional three-web
with a nonvanishing covector $a$ and with the integrable
transversal $a$-distribution $\Delta$ defined by this covector.
 Two-dimensional three-webs $W (3, 2, 1)$ cut
on integral surfaces $V^2$  of the transversal $a$-distribution $\Delta$
 by the foliations of the web $W (3, 2, 2)$ are hexagonal if and only if
the $a$-distribution $\Delta$ is one of four principal transversal
distributions of the pseudoconformal structure $CO (2, 2)$
associated with the web $W (3, 2, 2)$.
\end{theorem}

{\sf Proof.} {\em Sufficiency}. Using the same considerations
which we used in the proof of Theorem 8, we find that the
integrability of $\Delta$ and the hexagonality of $W (3, 2, 1)$
lead to conditions (52). Equations (52) and (56) give $b = 0$.
By (57), the last condition means that the transversal bivectors
 of the $a$-distribution $\Delta$ are principal.

 {\em Necessity}. If the $a$-distribution $\Delta$ is integrable and all
 bivectors $\Delta$ of the pseudoconformal structure $CO (2, 2)$
defined by the web $W (3, 2, 2)$ on $M^4$ are principal, then
we have the first equation of (53), $b^{1'}_{2' 2' 2'} = 0$,
and $b = b^i a_i = 0$. These two conditions imply that
$$
K = b^{2'}_{2' 2' 2'} = 0,
$$
i.e., the three-webs $W (3, 2, 1)$ cut
on integral surfaces $V^2$  of the transversal $a$-distribution $\Delta$
 by the foliations of the web $W (3, 2, 2)$ are hexagonal.
 \rule{3mm}{3mm}

\noindent {\em Authors' addresses}:\\

\noindent
\begin{tabular}{ll}
M. A. Akivis &V. V. Goldberg\\
Department of Mathematics &Department of Mathematical Sciences\\
Jerusalem College of Technology---Mahon Lev &  New
Jersey Institute of Technology \\
Havaad Haleumi St., P. O. B. 16031 & University Heights \\
 Jerusalem 91160, Israel &  Newark, N.J. 07102, U.S.A. \\
 & \\
 E-mail address: akivis@avoda.jct.ac.il & E-mail address:
 vlgold@m.njit.edu
 \end{tabular}

\end{document}